\documentclass{my_aims}

\usepackage{txfonts}
\def\typeofarticle{Research article}

\def\ppages{xxx--xxx}
\def\DOI{}
\def\Received{8 September 2021}
\def\Accepted{3 November 2021}
\def\Published{}

\usepackage{amsmath,amssymb,latexsym}
\usepackage{subfig}
\usepackage{psfrag}
\usepackage{enumitem}
\usepackage{multirow}

\newcommand{\diff}{\mathrm{d}}

\def\R{\mathbb{R}}
\def\Z{\mathbb{Z}}

\newtheorem{remark}{Remark}

\numberwithin{equation}{section}

\begin{document}
	

\title{Model-free based control of a HIV/AIDS prevention model}

\author{%
  Lo\"ic Michel\affil{1,}\affil{2}\corrauth,
  Cristiana J. Silva\affil{3} 
  and
  Delfim F. M. Torres\affil{3} 
}

\shortauthors{L. Michel, C. J. Silva, D. F. M. Torres}


\address{%
  \addr{\affilnum{1}}{\'Ecole Centrale de Nantes-LS2N, UMR 6004 CNRS, 44300 Nantes, France}
  \addr{\affilnum{2}}{Univ Lyon, INSA Lyon, Universit\'e Claude Bernard Lyon 1, 
  	                  \'Ecole Centrale de Lyon, CNRS, Amp\`ere, UMR 5005, 69621 Villeurbanne, France}
  \addr{\affilnum{3}}{Center for Research and Development in Mathematics and Applications (CIDMA),\\
		              Department of Mathematics, University of Aveiro, 3810-193 Aveiro, Portugal}}

\corraddr{loic.michel@ec-nantes.fr}


\begin{abstract}
Controlling an epidemiological model is often performed using optimal control 
theory techniques for which the solution depends on the equations 
of the control system, objective functional and possible state and/or control constraints. 
In this paper, we propose a model-free control approach based on an algorithm 
that operates in 'real-time' and drives the state solution according to a direct 
feedback on the state solution that is aimed to be minimized, and without knowing 
explicitly the equations of the control system. We consider a concrete epidemic 
problem of minimizing the number of HIV infected individuals, through 
the preventive measure \emph{pre-exposure prophylaxis (PrEP)} given to susceptible individuals. 
The solutions must satisfy control and mixed state-control constraints that represent 
the limitations on PrEP implementation. Our model-free based control algorithm allows 
to close the loop between the number of infected individuals with HIV and the supply 
of PrEP medication 'in real time', in such manner that the number of infected individuals 
is asymptotically reduced and the number of individuals under PrEP medication is below 
a fixed constant value. We prove the efficiency of our approach and compare 
the model-free control solutions with the ones obtained using a classical optimal 
control approach via Pontryagin maximum principle. The performed numerical simulations 
allow us to conclude that the model-free based control strategy highlights 
new and interesting performances compared with the classical optimal control approach.
\end{abstract}

\keywords{HIV/AIDS model; model-free control; real-time control; intelligent PID control; tracking.}

\maketitle


\section{Introduction}

The SICA (Susceptible--Infected--Chronic--AIDS) compartmental model for HIV/AIDS transmission
dynamics was proposed by Silva and Torres in their seminal paper \cite{SilvaTorres:TBHIV:2015}.
Since then, the model has been extended to stochastic systems of differential equations 
\cite{D:S:T:2018,MR4261252}, to fractional-order \cite{S:T:2019,MR4279814} as well as discrete-time 
dynamics \cite{V:T:2021}, and applied with success to describe very different HIV/AIDS epidemics, 
like the ones in Cape Verde \cite{SilvaTorres:EcoComplexity,SilvaTorres:PrEP:DCDS:2018} 
or Morocco \cite{L:M:M:S:T:Y:2019}. For a survey see \cite{MyID:455}. 

In this work, we consider a SICA epidemic problem of controlling the transmission 
of the human immunodeficiency viruses (HIV), by considering not only the medical treatment 
with multiple antiretroviral (ART) drugs, but also the pre-exposure prophylaxis (PrEP), 
which are medicines taken to prevent getting HIV infection. According to the
\emph{Centers for Disease Control and Prevention}, PrEP is highly effective for preventing HIV, 
when taken as prescribed, and reduces the risk of getting HIV from sex and from injection 
drug use by about 99\% and 74\%, respectively \cite{cdc}. 

The main contribution of this work is to propose a model-free based control algorithm 
that closes the loop between the infected individuals with HIV and PreP medication, 
in such manner that the medication is driven in 'real-time', according to the number 
of infected individuals that has to be asymptotically reduced. We highlight 
that {\it model-free control} offers the advantages of a simple  
Proportional--Integral--Derivative (PID) controller 
in the framework of model free design, that is, one whose parameters that can be easily 
tuned without a precise knowledge of the controlled epidemiological model.

The model-free control methodology, originally proposed 
by Fliess and Join in \cite{Fliess}, has been designed 
to control {\it a priori} any ``unknown'' dynamical system 
in a ``robust'' manner, and is referred to as ``a self-tuning regulator'' 
in \cite{Astrom}. This control law can be considered as an alternative 
to PI and PID controllers \cite{FliessJoin_2021} and the performances 
are really satisfactory taking into account that the control is calculated 
based only on the information provided by the controlled input and the measured 
output signal of the controlled systems. This control law has been extensively 
and successfully applied to control many nonlinear processes: see, e.g., 
\cite{Bara,Fliess,Hamiche} and the references therein. In particular, 
some applications have been dedicated to the control of chemistry 
and biological processes \cite{Bara2,Bara,Bernier,Ridha,Tebbani}, 
including the development of an artificial pancreas \cite{glycemia}.
A derivative-free-based version of this control algorithm has been 
proposed by the first author in \cite{michel2018}, for which some 
interesting capabilities of online optimization have been highlighted. 
To the best of our knowledge, the application of model-free control 
to SICA modeling has never been discussed before.

Here, we compare the solutions obtained by the model-free control method 
with the corresponding solutions of an optimal control problem 
for HIV/AIDS transmission from \cite{SilvaTorres:PrEP:DCDS:2018}, 
which has a mixed state-control constraint. In \cite{SilvaTorres:PrEP:DCDS:2018}, 
the control system is based on a SICAE (Susceptible, HIV-Infected, 
Chronic HIV-infected under ART, AIDS-symptomatic individuals, 
E -- under PrEP medication) model for the transmission of HIV 
in a homogeneously mixing population. The control $u$ represents the fraction of susceptible 
individuals under PrEP, with $0 \leq u(t) \leq 1$, that is, when $u(t)=0$, no susceptible individual 
takes PrEP at time $t$, and when $u(t) = 1$ all susceptible individuals are taking PrEP 
at time $t$. The mixed state-control constraint refers to the fact that only people 
who are HIV-negative and at a very high risk of HIV infection should take PrEP, 
and also to the high costs of PrEP medication. Therefore, the number of susceptible 
individuals that takes PrEP, at each day, must be bounded by a positive constant. 
Moreover, the cost functional, which is aimed to be minimized, represents a balance 
between the number of HIV infected individuals and the costs associated with PrEP implementation. 

The paper is structured as follows. In Section~\ref{sec:materials}, we propose a model-free 
control method and the procedure to minimize the HIV infected cases is described. 
In Section~\ref{sec:results}, we present some numerical simulations and provide 
a comparison of the results obtained using the model-free based approach with the ones 
in \cite{SilvaTorres:PrEP:DCDS:2018} from the Pontryagin maximum principle. 
Section~\ref{sec:discussion} discusses and compares the results. 
Finally, some concluding remarks and possible directions
for future work are given in Section~\ref{sec:conclusion}. 


\section{Materials and methods}
\label{sec:materials}

In this section, we propose our model-free based control method 
and apply it to an epidemiological problem of minimizing 
HIV-infected individuals. 


\subsection{Principle of the model-free based control}

Model-free based control was introduced in 2008 and 2009 
by Fliess and Join in \cite{Fliess2008a,Fliess2009}. It is an alternative 
technique to control complex systems based on elementary continuously 
updated local modeling via unique knowledge of the input-output behavior. 
The key feature of this approach lies in the fact that the control system, 
which might be highly nonlinear and/or time-varying, is taken into account 
without any modeling procedure \cite{Nascimento}. The model-free based 
control approach has been successfully implemented in concrete applications 
to diverse fields, ranging from intelligent transportation systems 
to energy management, etc., see \cite{Fliess} and references cited therein. 
To the best of our knowledge, no one as yet used this approach 
in the context of epidemiology.

Consider a nonlinear dynamical system $f : u \mapsto y$ to control:
\begin{equation}
\label{eq:gen_sys}
\left\{ \begin{array}{l}
\dot{ \boldsymbol{x} } = f( \boldsymbol{x},u), \\
y = g( \boldsymbol{x}), 
\end{array} \right.
\end{equation}
where $f$ is the function describing the behavior of a nonlinear system and 
$\boldsymbol{x} \in \R$ is the state vector. The para-model control is an 
application $\mathcal{C}_{\pi} : (y, y^*) \mapsto u$, whose purpose 
is to control the output $y$ of~\eqref{eq:gen_sys} 
following an output reference $y^*$. In simulation, the system~\eqref{eq:gen_sys} 
is controlled in its ``original formulation'', without any modification or linearization.

For any discrete moment $t_k, \, k \in \mathbb{N}^*$, one defines the discrete controller 
$\mathcal{C}_{\pi} : (y, y^*) \mapsto u$ as an integrator associated to a numerical series 
$(\Psi_k)_{k \in \Z}$, symbolically represented by
\begin{equation} 
\label{eq:iPI_discret_nm_eq}
u_k = \mathcal{C}_{\pi} (y_{k}, y^*_{k}) = \Psi_k 
\cdot \int_0^t K_i (y^\ast_{k} - y_{k-1}) \, d \, \tau
\end{equation}
with the recursive term
\begin{equation*}
\Psi_k = \Psi_{k-1} + {K_p} ( k_\alpha e^{-k_\beta k} - y_{k-1}),
\end{equation*}
where $y^\ast$ is the output (or tracking) reference trajectory; 
$K_p$ and $K_i$ are real positive tuning gains; $\varepsilon_{k-1} 
= y^\ast_{k} - y_{k-1}$ is the tracking error; and $k_\alpha e^{-k_\beta k}$ 
is an initialization function, where $k_\alpha$ and $k_\beta$ are real positive constants. 
In practice, the integral part is discretized using, for example, Riemann sums. 

The set of the $\mathcal{C}_{\pi}$-parameters of the controller, 
is defined as the set of the tuning coefficients $\{K_p, K_i, k_\alpha, k_\beta\}$.

\subsection{Methodology}

Consider the problem of minimizing the number of infected individuals with HIV, 
given by the state trajectory $I$, through the control measure $u$ associated 
to PrEP medication, and satisfying a mixed state-control constraint 
$S u \leq \eta$, where $\eta$ is a positive constant. 

\medskip

The control sequence is divided into two steps, in order to manage, separately, 
the increasing $u$ transient and the associated decreasing $u$ transient 
that must decrease afterwards both the infected cases and the control 
input $u$ to lower values on $u$ and $I$. Concerning the control input $u$, 
it is expected that $u(t > t_f) = 0$, where $t_f$ 
is the time from which the medication is stopped. 

\begin{itemize}
\item A \emph{first sequence}, associated to the setting up of the medication, 
aims to progressively increase the medication (and thus start decreasing the infected states) 
until a certain threshold of infected cases is reached, above which the number of infected 
cases could be considered as stabilizable around a low value. This control sequence can 
be managed by a simple $L$-linear or $Q$-quadratic slope, such as
\begin{equation}
\mathrm{(slope)}  \quad  u(t) = u_0 + L \cdot t  \qquad  
\qquad  \mathrm{(quadratic)} \quad   u(t) = u_0 + Q \cdot t^2
\end{equation}
driven in open-loop, {\it i.e.}, independently from any feedback of the infected state, 
that increases gradually the medication while satisfying the constraint on $S u \leq \eta$ 
due to the remaining high level of the susceptible cases. It appears crucial to accelerate 
the medication at the beginning, in order to reach rapidly $\gamma_{lim}$ and allow a strong 
decreasing of the infected cases. This point will be discussed later in Section~\ref{sec:results}.
    
\item Denote $u_{max} = \max u(t)$, the maximum value of $u$. The \emph{second sequence} 
is associated to the decrease of the medication until the infected state is stabilized 
around a low value. This sequence is managed by our proposed model-free based control 
that interacts, in real-time, with the number of infected cases and, consequently, 
calculates the ``optimal'' medication $u$ in order to decrease and stabilize 
the infected state. According to \eqref{eq:iPI_discret_nm_eq}, the control reads:
\begin{equation*} 
u_k = \Psi_k \, \cdot \int_0^t K_i (I^\ast_{k} - I_{k-1}) \, d \, \tau  
\qquad \mathrm{with}  \qquad \Psi_k 
= \Psi_{k-1} + {K_p} ( k_\alpha e^{-k_\beta k} - I_{k-1}),
\end{equation*}
where $I^\ast$ denotes the infected cases reference that practically 
can be chosen as $I^\ast = \min I(t)$, where the minimum value of $I$ can be reached online, 
updating the tracking reference and, therefore, ensuring that the control law is ``synchronized'' 
on the lowest value that can be reachable. It is worth to note that, depending on the behavior 
of the closed-loop, a saturation is added to bound the controlled $u$: $0 \leq u \leq  u_{max} \leq 1$.

\end{itemize}

\begin{figure}[!h]
\centering
\includegraphics[width=14cm]{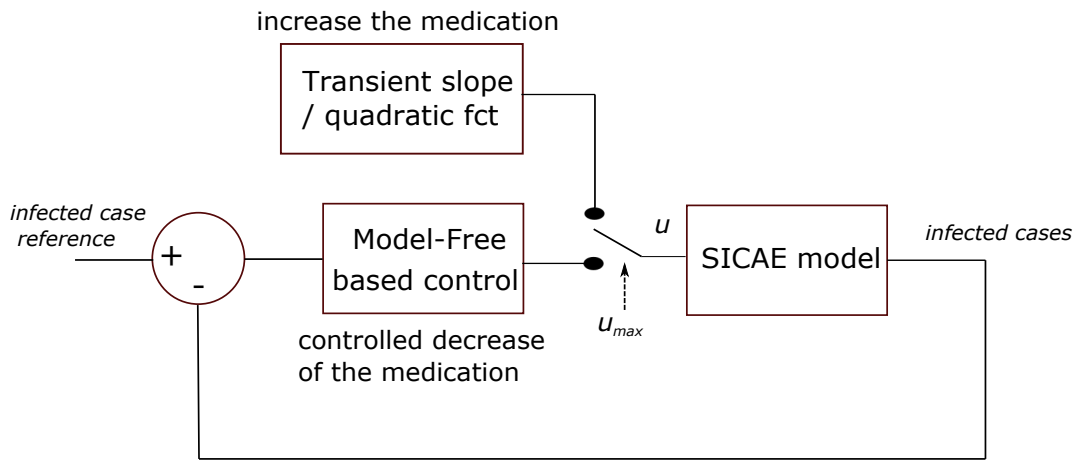}
\caption{Proposed scheme to control a nonlinear system.}
\label{fig:CSM_gen}
\end{figure}

The implementation of the control scheme is depicted in Figure~\ref{fig:CSM_gen}, 
where the control sequence starts from the transient slope or the quadratic function 
and then, once $u_{max}$ is reached, switches to the proposed model-free based controller.
The parameters of the control sequence to be adjusted, comprise the 
$\mathcal{C}_{\pi}$-parameters set of the model-free control algorithm 
and the $(u_0, L, Q)$ parameters of the first sequence, 
depending if a linear or a quadratic slope is involved.

\paragraph{Additional constraints.} 

The slope can be adjusted according to a state-control constraint that determines the maximum number 
of susceptible individuals that take PrEP medication, at each instant of time. This constraint reads 
as $S(t) \cdot u(t) < \gamma_{max}$, for all $t$, where $\gamma_{max}$ is the corresponding upper bound.

To properly tune each sequence, in order to satisfy both the state-control constraint 
as well as to minimize the cost criteria, a derivative-free based optimization 
procedure can be applied \cite{Porcelli}.

\subsection{Classical optimal control problem applied to the HIV/AIDS SICAE model with a mixed state-control constraint}
\label{subsec:2.3}

Consider the mathematical model for HIV/AIDS transmission in a homogeneously mixing population 
proposed in \cite{SilvaTorres:PrEP:DCDS:2018} and based in 
\cite{SilvaTorres:TBHIV:2015,SilvaTorres:EcoComplexity}.

The model subdivides human population into five mutually-exclusive 
compartments: susceptible individuals ($S$); 
HIV-infected individuals with no clinical symptoms of AIDS 
(the virus is living or developing in the individuals 
but without producing symptoms or only mild ones) 
but able to transmit HIV to other individuals ($I$); 
HIV-infected individuals under ART treatment (the so called 
chronic stage) with a viral load remaining low ($C$); 
HIV-infected individuals with AIDS clinical symptoms ($A$); 
individuals that are under PrEP medication ($E$). 
The total population at time $t$, denoted by $N(t)$, is given by
$N(t) = S(t) + I(t) + C(t) + A(t) + E(t)$.
The model assumptions are the following \cite{SilvaTorres:TBHIV:2015,SilvaTorres:PrEP:DCDS:2018}. 
Effective contact with people infected with HIV is at a rate $\lambda$, given by
\begin{equation*}
\lambda = \frac{\beta}{N} \left( I + \eta_C \, C  + \eta_A  A \right),
\end{equation*}
where $\beta$ is the effective contact rate for HIV transmission.
The modification parameter $\eta_A \geq 1$ accounts for the relative
infectiousness of individuals with AIDS symptoms, in comparison to those
infected with HIV with no AIDS symptoms. Individuals with AIDS symptoms
are more infectious than HIV-infected individuals (pre-AIDS) because
they have a higher viral load and there is a positive correlation
between viral load and infectiousness \cite{art:viral:load}. 
On the other hand, $\eta_C \leq 1$ translates the partial restoration 
of immune function of individuals with HIV infection 
that use ART correctly \cite{AIDS:chronic:Lancet:2013}.
All individuals suffer from natural death, at a constant rate $\mu$. 
HIV-infected individuals, with and without AIDS symptoms, have access to ART treatment. 
HIV-infected individuals with no AIDS symptoms $I$ progress to the class 
of individuals with HIV infection under ART treatment $C$ at a rate $\phi$, 
and HIV-infected individuals with AIDS symptoms are treated for HIV, at rate $\alpha$.
HIV-infected individuals with AIDS symptoms $A$, 
that start treatment, move to the class of HIV-infected individuals $I$, 
moving to the chronic class $C$ only if the treatment is maintained. 
HIV-infected individuals with no AIDS symptoms $I$ that do not take 
ART treatment progress to the AIDS class $A$, at rate $\rho$. 
Individuals in the class $C$ that stop ART medication are transferred 
to the class $I$, at a rate $\omega$. 
Only HIV-infected individuals with AIDS symptoms $A$ 
suffer from an AIDS induced death, at a rate $d$. 
The proportion of susceptible individuals that takes PrEP is denoted by $\psi$.
It is assume that PrEP is effective, so that all susceptible individuals 
under PrEP treatment are transferred to class $E$. The individuals 
that stop PrEP become susceptible individuals again, at a rate $\theta$. 
Susceptible individuals are increased by the recruitment rate $\Lambda$.  

Such model is given by the following system of ordinary differential equations:
\begin{equation}
\label{eq:model:PreP}
\begin{cases}
\dot{S}(t) = \Lambda - \frac{\beta \left( I(t) + \eta_C \, C(t)  
+ \eta_A  A(t) \right)}{N(t)} S(t) - \mu S(t) - \psi S(t) + \theta E(t),\\[0.2 cm]
\dot{I}(t) = \frac{\beta \left( I(t) + \eta_C \, C(t)  
+ \eta_A  A(t) \right)}{N(t)} S(t) - (\rho + \phi + \mu) I(t) 
+ \alpha A(t)  + \omega C(t), \\[0.2 cm]
\dot{C}(t) = \phi I(t) - (\omega + \mu)C(t),\\[0.2 cm]
\dot{A}(t) =  \rho \, I(t) - (\alpha + \mu + d) A(t) ,\\[0.2 cm]
\dot{E}(t) = \psi S(t) - (\mu + \theta) E(t).
\end{cases}
\end{equation}

Recall that PrEP medication should only be administrated to people who are HIV-negative 
and at very high risk for HIV infection. Moreover, PrEP is highly expensive 
and it is still not approved in many countries.  Therefore, 
the number of individuals that should take PrEP should be limited at each instant 
of time for a fixed interval of time $[0, t_f]$ \cite{SilvaTorres:PrEP:DCDS:2018}.
The optimal control problem proposed in \cite{SilvaTorres:PrEP:DCDS:2018}, 
and considered in this paper for comparison of results, 
takes into account this health public problem. 

The main goal of the optimal control problem is to determine 
the PrEP strategy $\psi$ that minimizes the number of individuals with pre-AIDS 
HIV-infection $I$ as well as the costs associated with PrEP. Let the fraction 
of individuals that takes PrEP, at each instant of time, be a control function, 
that is, $\psi \equiv u(t)$ with $t \in [0, t_f]$, and assume that the total population 
$N$ is constant: the recruitment rate is proportional to the natural death rate, 
$\Lambda = \mu N$, and there are no AIDS-induced deaths ($d=0$). 
The controlled model is given by
\begin{equation}
\label{eq:model:PreP:control}
\begin{cases}
\dot{S}(t) = \mu N - \frac{\beta}{N} \left( I(t) + \eta_C \, C(t)  
+ \eta_A  A(t) \right) S(t) - \mu S(t) - S(t) u(t) + \theta E(t),\\[0.2 cm]
\dot{I}(t) = \frac{\beta}{N} \left( I(t) + \eta_C \, C(t)  
+ \eta_A  A(t) \right) S(t) - (\rho + \phi + \mu)I(t) + \alpha A(t)  + \omega C(t), \\[0.2 cm]
\dot{C}(t) = \phi I(t) - (\omega + \mu)C(t),\\[0.2 cm]
\dot{A}(t) =  \rho \, I(t) - (\alpha + \mu) A(t) ,\\[0.2 cm]
\dot{E}(t) = S(t) u(t) - (\mu + \theta) E(t) \, .
\end{cases}
\end{equation}

\begin{remark} 
All the parameters of the SICAE model \eqref{eq:model:PreP:control} are fixed with the exception 
of the control function $u(t)$. This system is deterministic and there is no uncertainty. However, 
the model-free based control proposed in this paper does not use these equations. They are only 
needed in the classical optimal control approach that is used here for comparison. 
The sensitivity analysis of the parameters of the SICA model, which is in the basis 
of the SICAE model \eqref{eq:model:PreP}, was studied before 
in \cite{SilvaTorres:EcoComplexity}. 
\end{remark}

The classical optimal control problem proposed in \cite{SilvaTorres:PrEP:DCDS:2018}, 
and that is considered here in comparison with the model-free control method, 
considers the cost functional 
\begin{equation}
\label{eq:cost}
J(u) = \int_0^{t_f} \left[ w_1 I(t) + w_2 u^2(t) \right]  \, dt,
\end{equation}
where the constants $w_1$ and $w_2$ represent the weights associated with 
the number of HIV infected individuals $I$ and on the cost associated with 
the PrEP prevention treatment, respectively.
It is assumed that the control 
function $u$ takes values between 0 and 1. When $u(t) = 0$, no susceptible 
individual takes PrEP at time $t$; if $u(t)=1$, then all susceptible individuals 
are taking PrEP at time $t$. Let $\gamma_{max}$ denote the total number of susceptible 
individuals under PrEP for a fixed time interval $[0, t_f]$. This constraint is represented by
\begin{equation}
\label{eq:constraint}
S(t) u(t) \leq \gamma_{max} \, , \quad \gamma_{max} \geq 0 \, , \, \,
\text{for almost all} \, \, t \in [0, t_f] \, ,
\end{equation}
which should be satisfied at almost every instant of time during the whole PrEP program.

Let 
\begin{equation*}
x(t) =(x_1(t), \ldots,  x_5(t))
=\left( S(t), I(t), C(t), A(t), E(t) \right) 
\in {\mathbb{R}}^5.
\end{equation*}
The classical optimal control problem proposed in \cite{SilvaTorres:PrEP:DCDS:2018} 
consists to find the optimal trajectory $\tilde{x}$, 
associated with the control $\tilde{u}$, satisfying the control system 
\eqref{eq:model:PreP:control}, the initial conditions
\begin{equation*}
x(0) = (x_{10}, x_{20}, x_{30}, x_{40}, x_{50}), 
\quad \text{with} \quad x_{10} \geq 0,\,  x_{20} \geq 0, \, x_{30} \geq 0, \,  
x_{40} \geq 0, \,  x_{50} \geq 0,
\end{equation*}
the constraint \eqref{eq:constraint}, and where the control 
$\tilde{u} \in \Omega$,
\begin{equation}
\label{eq:admiss:control}
\Omega = \biggl\{ u(\cdot) \in L^{\infty}(0, t_f) \,
| \,  0 \leq u (t) \leq 1  \biggr\},
\end{equation}
minimizes the objective functional \eqref{eq:cost}.


\subsection{Comparative study and cost criteria definitions}

We evaluate the accuracy of our model-free proposed approach compared to the classical optimal 
control one when applied to the SICAE model described in Section~\ref{subsec:2.3}.
 
Let $T_e$ denote the final time of the $u$ treatment such as 
$u(t \geq T_e) = 0$ and $I(T_e)$ denotes the final value of the state ``Infected'' at the time $T_e$.
Regarding the ``energy'' of the control input $u$ with respect to the behavior 
of the infected cases and the period of time $T_e$ for which the medication $u$ 
is in effect, i.e., $u(t \geq T_e) = 0$, let us consider the following cost criteria: 
\begin{itemize}
\item \underline{Cost criterion for performances over the final time $T_e$:}
\begin{equation}
    J_{u+I} = \int_0^{T_e} u^2 + I^2 \, \diff \tau,
\end{equation}
\begin{equation}
    J_{I} = \int_0^{T_e} I^2 \, \diff \tau.
\end{equation}
    
\item  \underline{Time-pondered cost criterion:} 
\begin{equation}
J_{u+I}^{Te} = \int_0^{T_e} \tau \, ( u^2 + I^2 ) \, \diff \tau,
\end{equation}
which takes into account the effective period needed to stabilize 
the infected case, {\it i.e.}, the period for which $u > 0$.

\end{itemize}	


\section{Numerical simulations and results}
\label{sec:results}

To perform the numerical simulations, we consider the following parameter values, 
borrowed from \cite{SilvaTorres:PrEP:DCDS:2018}: $N = 10200$, $\mu = 1/69.54$, $\beta = 0.582$, 
$\eta_C = 0.04$, $\eta_A = 1.35$, $\theta = 0.001$,  $\omega = 0.09$, $\rho = 0.1$, 
$\phi = 1$ and $\alpha = 0.33$. The weight constants take the values $w_1 = w_2 = 1$.

The initial conditions are given by
\begin{equation*}
S(0) = 10000, \, \quad I(0) = 200, \, \quad C(0) = 0, \,  \quad 
A(0) = 0 \quad \text{and} \quad E(0) = 0,
\end{equation*}
and the mixed state-control constraint is
\begin{equation}
\label{eq:constraint:num}
S(t) \, u(t) \leq 2000 \, , \, \,
\text{for almost all} \, \, t \in [0, t_f].
\end{equation}

\begin{figure*}[!ht]
\centering
\subfloat[Evolution of the infected state $I$ versus time (in years).]{
\includegraphics[width=0.7\textwidth]{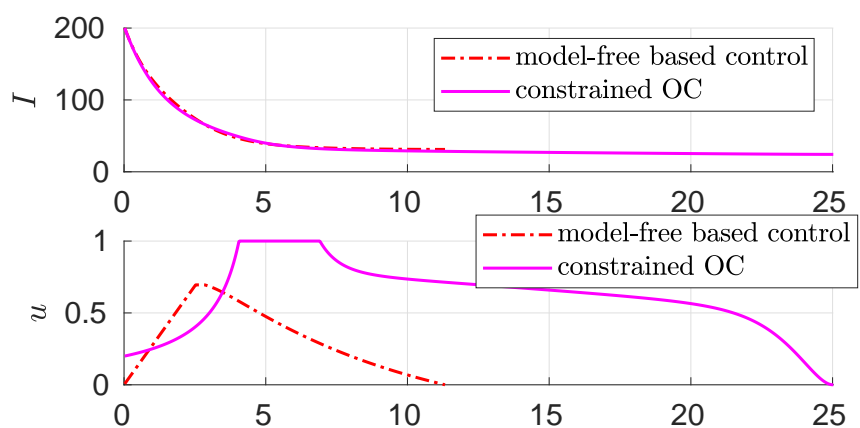}
\label{fg:non_opt:a}}
\hfil
\subfloat[Evolution of the controlled medication $u$ versus time (in years).]{
\includegraphics[width=0.7\textwidth]{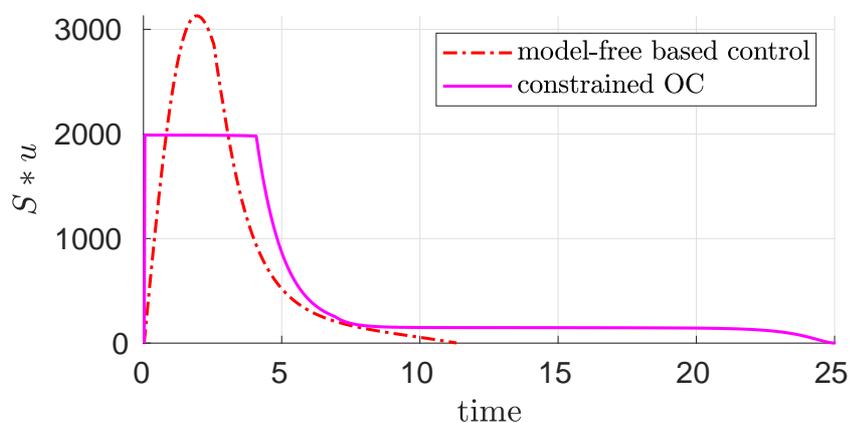}
\label{fg:non_opt:b}}
\caption{Evaluation of the unconstrained model-free based control 
with slope as first sequence.}
\label{fg:non_opt}
\end{figure*}

\begin{figure*}[!ht]
\centering
\subfloat[Evolution of the infected state $I$ versus time (in years).]{
\includegraphics[width=0.7\textwidth]{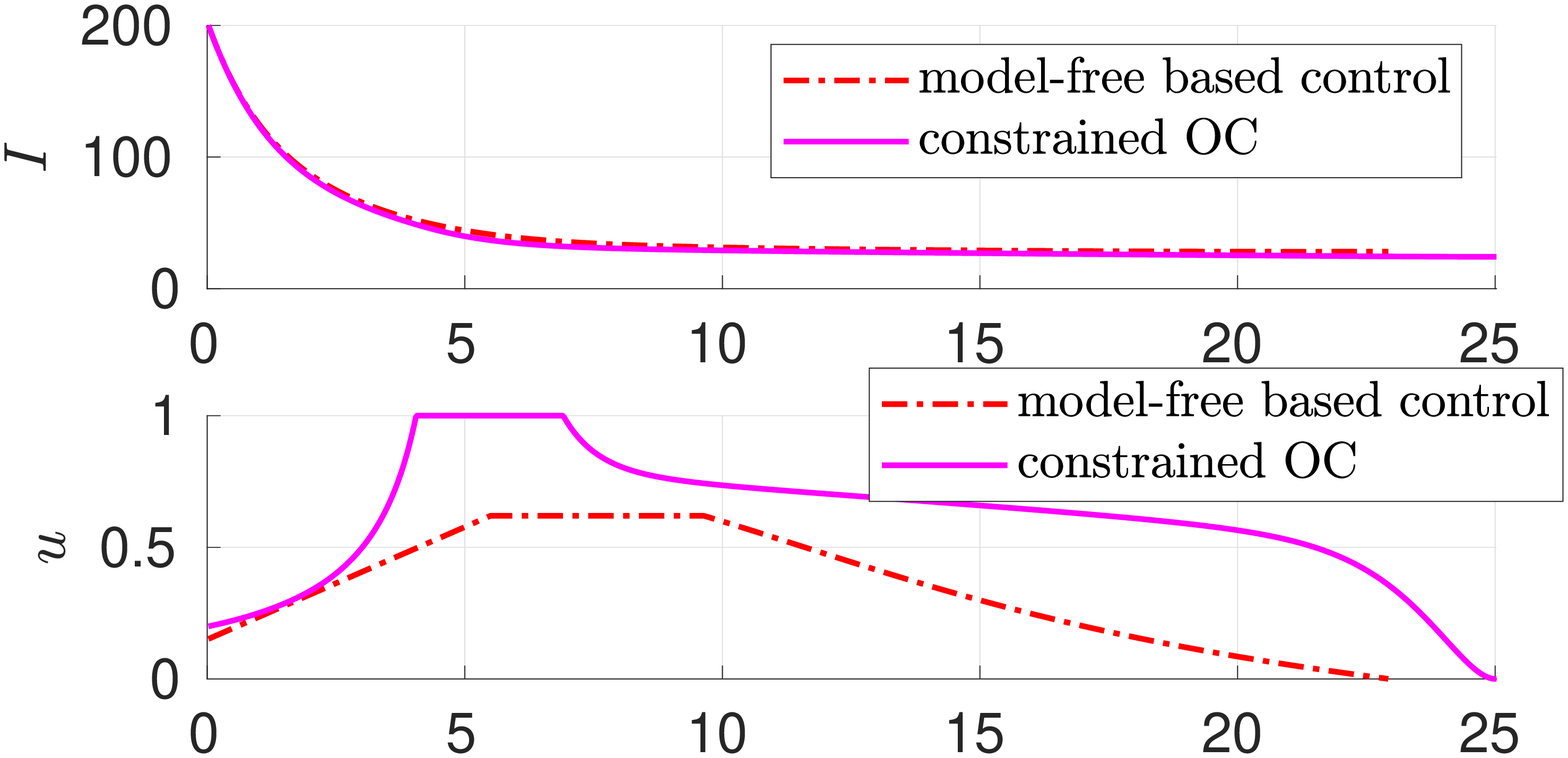}
\label{fg:opt_II:a}}
\hfil
\subfloat[Evolution of the controlled medication $u$ versus time (in years).]{
\includegraphics[width=0.7\textwidth]{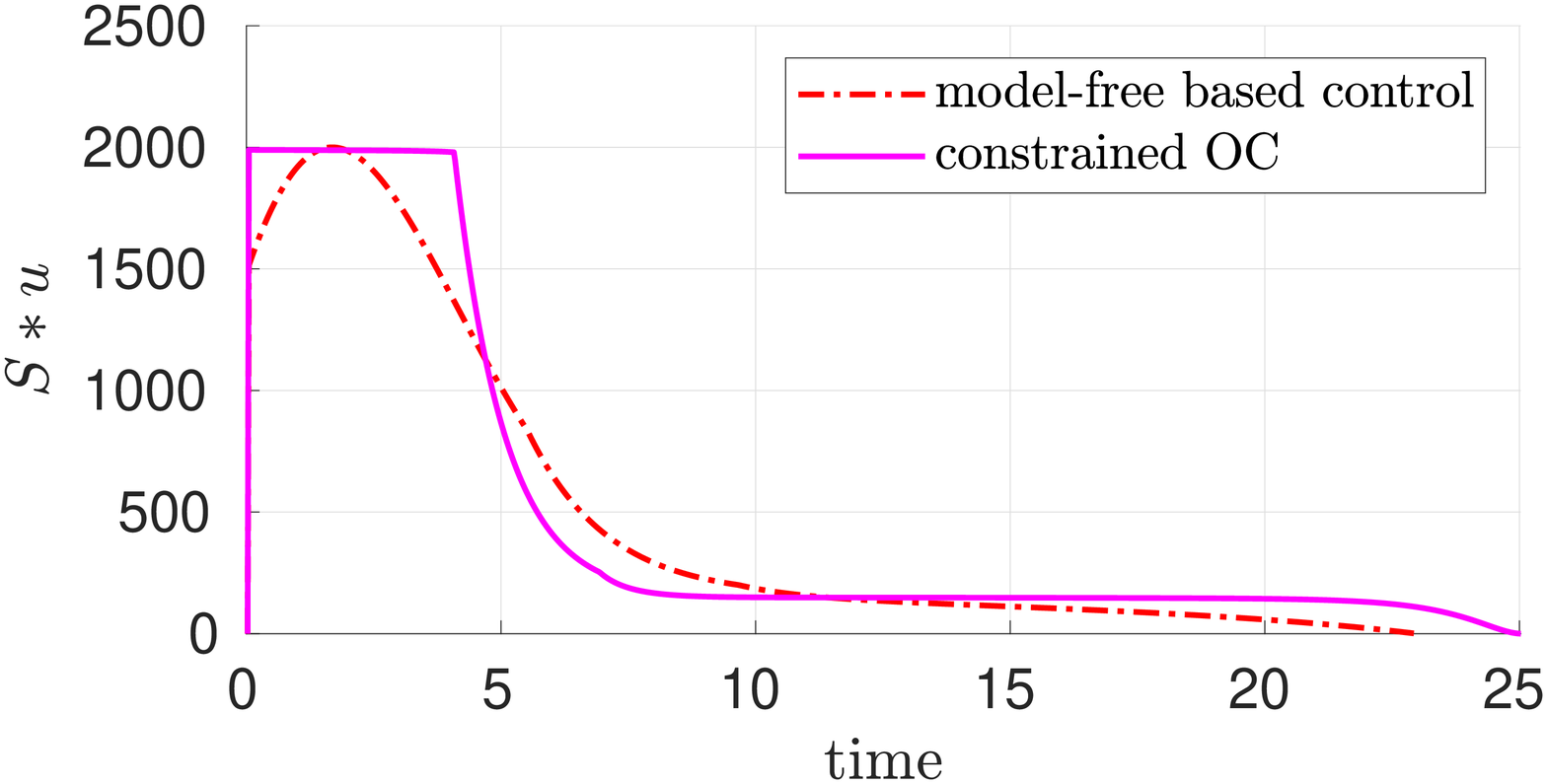}
\label{fg:opt_II:b}}
\caption{Evaluation of the constrained model-free based control: 
slope as first sequence including the constraint \eqref{eq:constraint:num} -- {\it case I}.}
\label{fg:opt_II}
\end{figure*}

\begin{figure*}[!ht]
\centering
\subfloat[Evolution of the infected state $I$ versus time (in years).]{
\includegraphics[width=0.7\textwidth]{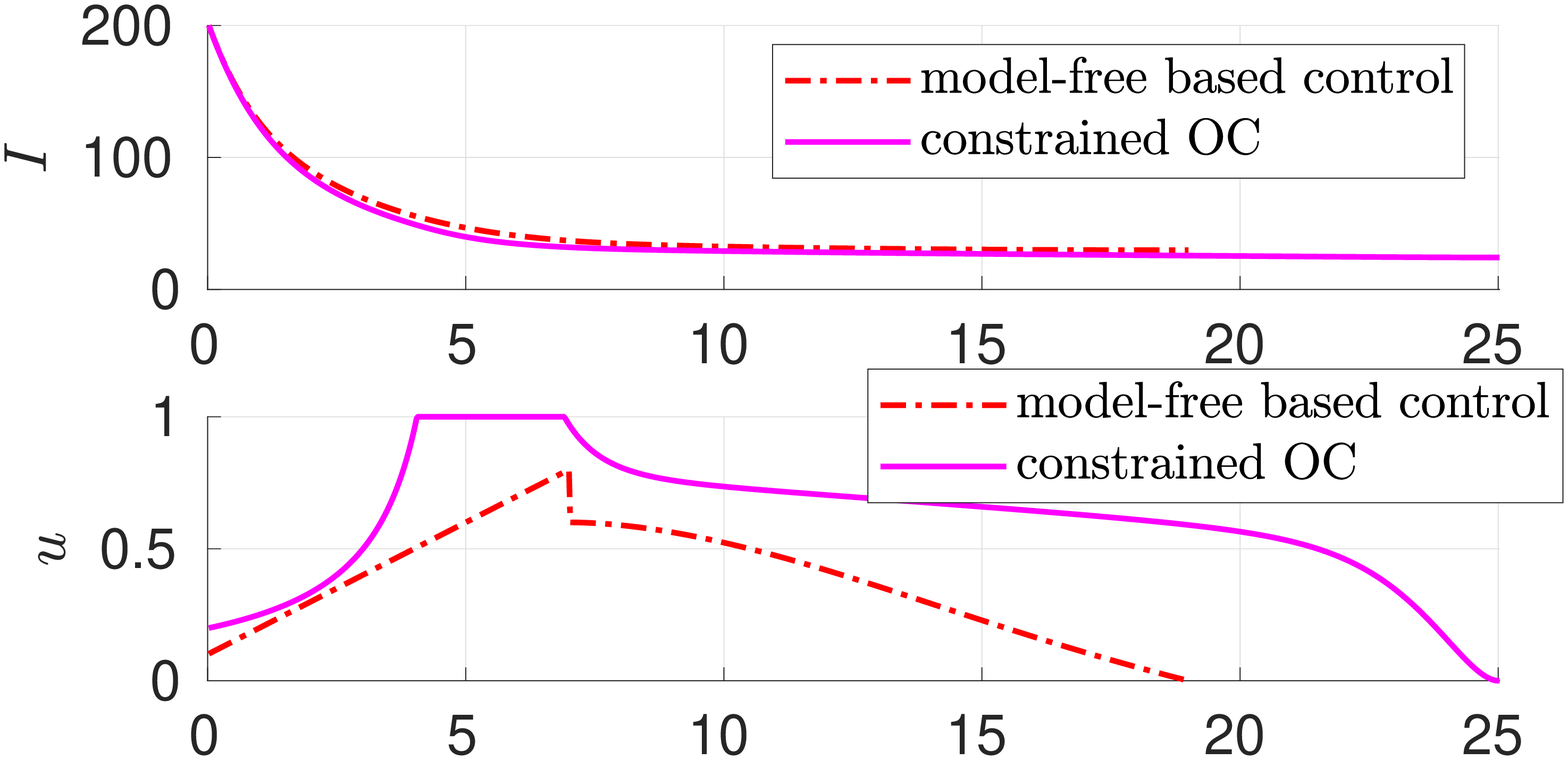}
\label{fg:opt_III:a}}
\hfil
\subfloat[Evolution of the controlled medication $u$ versus time (in years).]{
\includegraphics[width=0.7\textwidth]{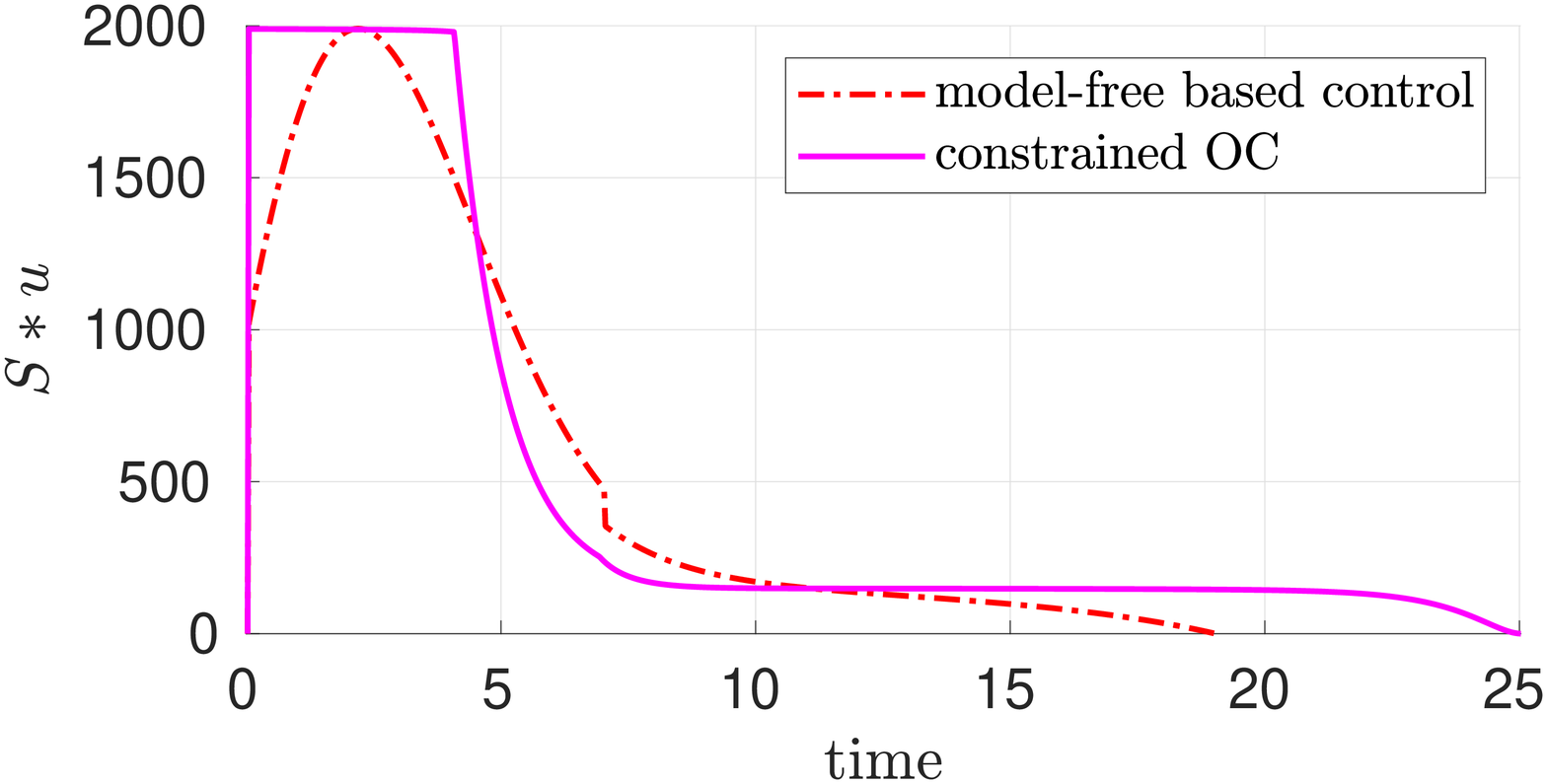}
\label{fg:opt_III:b}}
\caption{Evaluation of the constrained model-free based control: 
slope as first sequence including the constraint \eqref{eq:constraint:num} -- {\it case II}.}
\label{fg:opt_III}
\end{figure*}

\begin{figure*}[!ht]
\centering
\subfloat[Evolution of the infected state $I$ versus time (in years).]{
\includegraphics[width=0.7\textwidth]{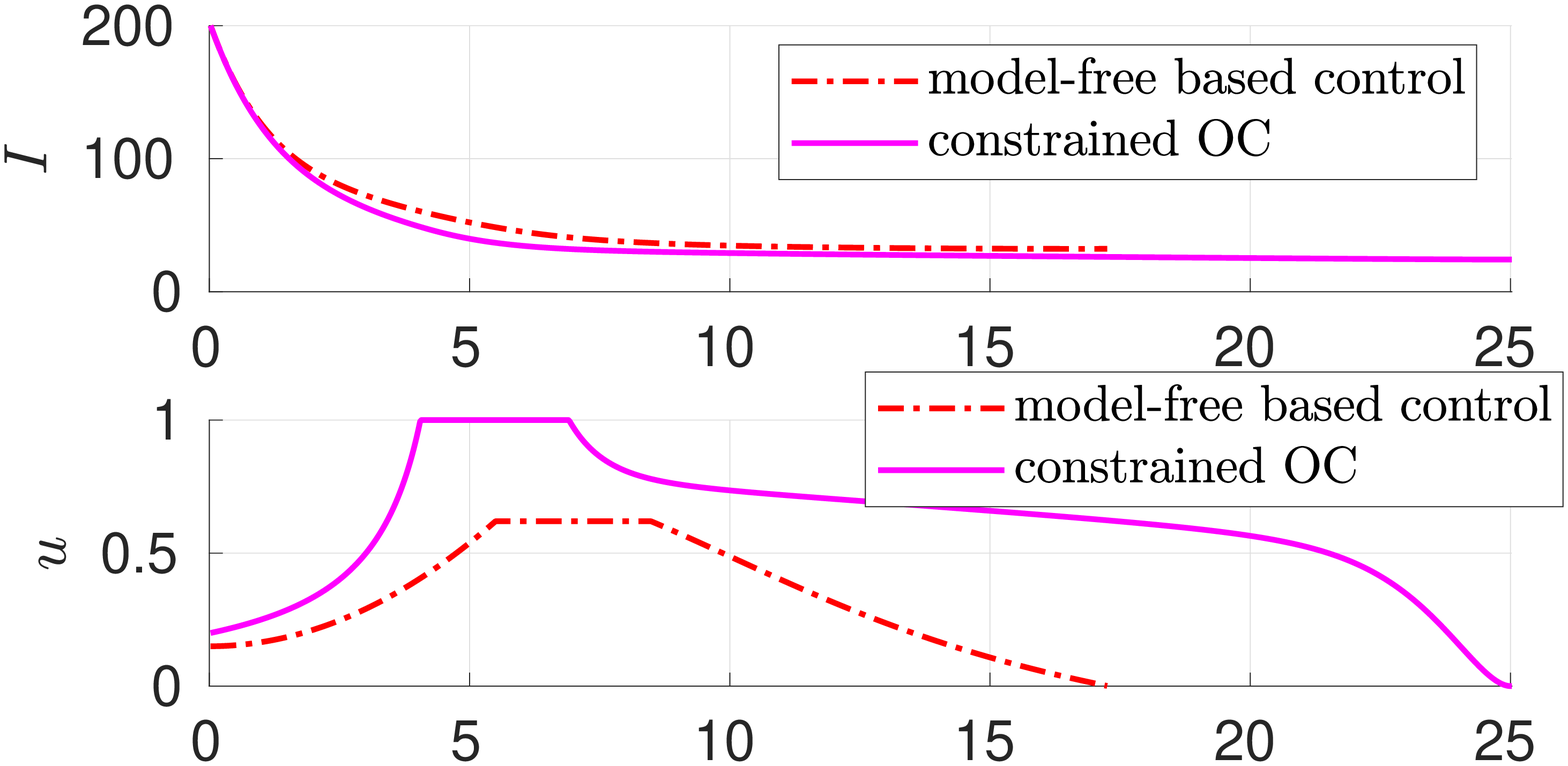}
\label{fg:opt_IV:a}}
\hfil
\subfloat[Evolution of the controlled medication $u$ versus time (in years).]{
\includegraphics[width=0.7\textwidth]{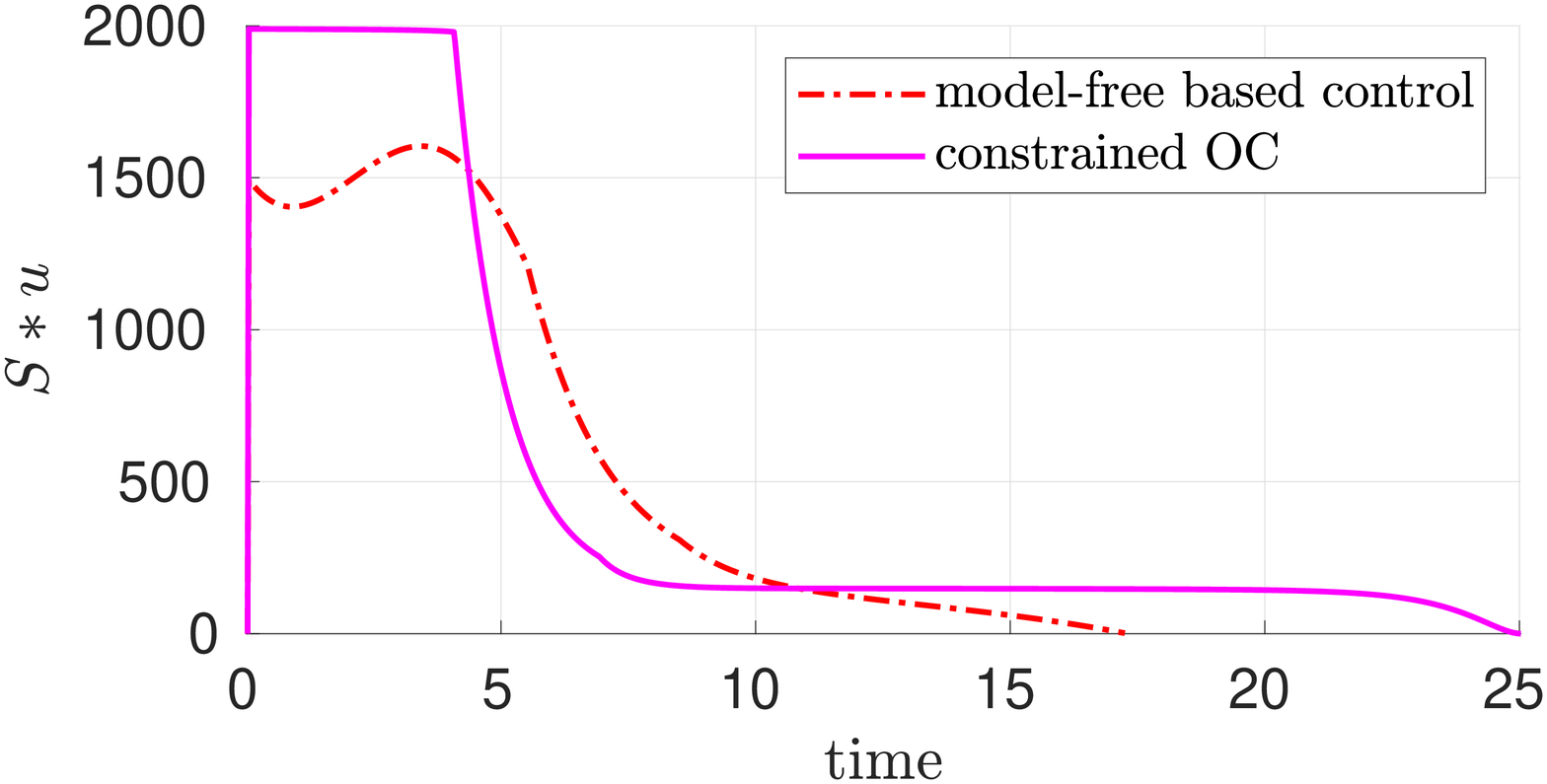}
\label{fg:opt_IV:b}}
\caption{Evaluation of the constrained model-free based control: 
quadratic as first sequence including the constraint \eqref{eq:constraint:num} -- {\it case I}.}
\label{fg:opt_IV}
\end{figure*}

\begin{figure*}[!ht]
\centering
\subfloat[Evolution of the infected state $I$ versus time (in years).]{
\includegraphics[width=0.7\textwidth]{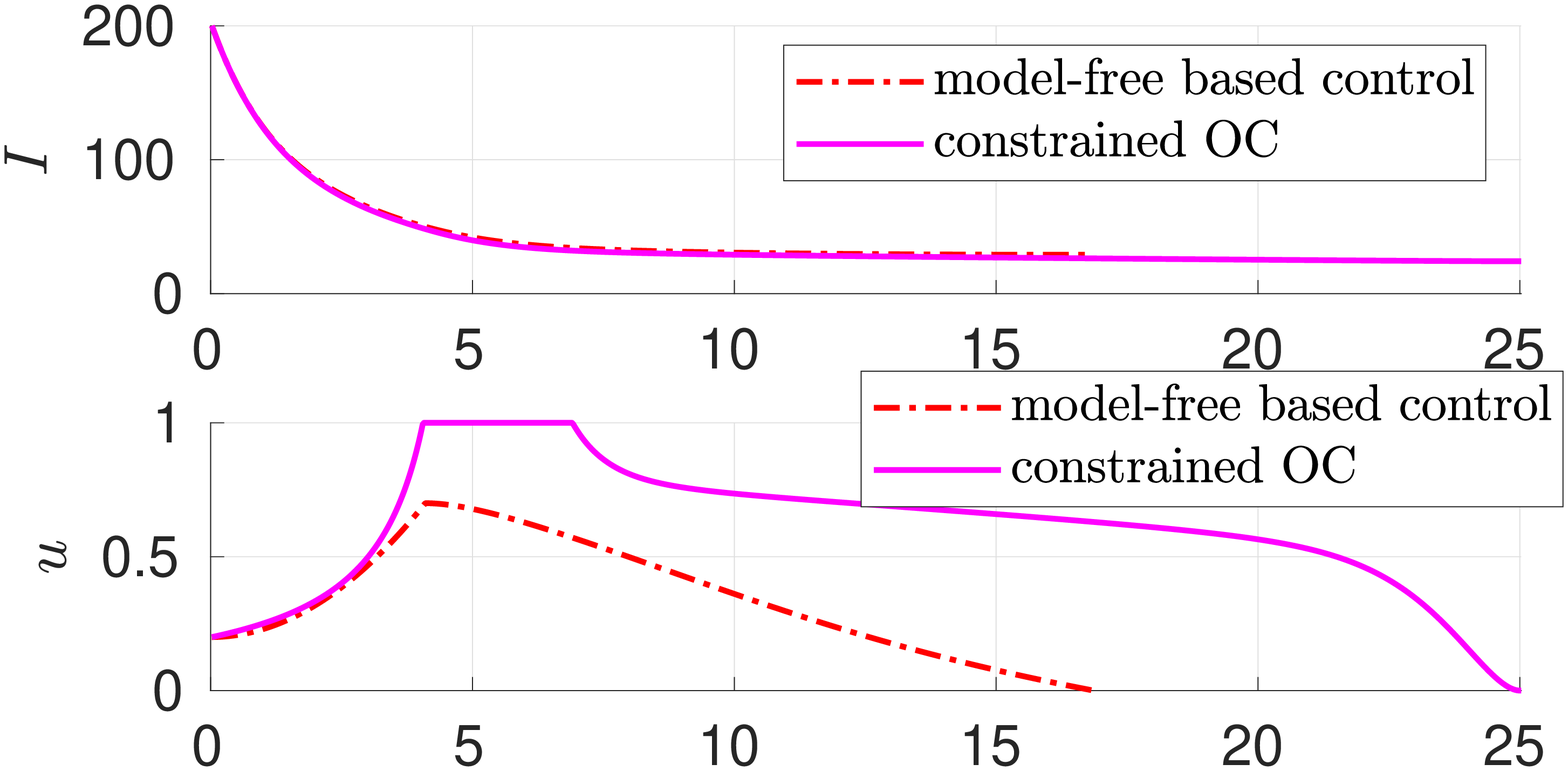}
\label{fg:opt_V:a}}
\hfil
\subfloat[Evolution of the controlled medication $u$ versus time (in years).]{
\includegraphics[width=0.7\textwidth]{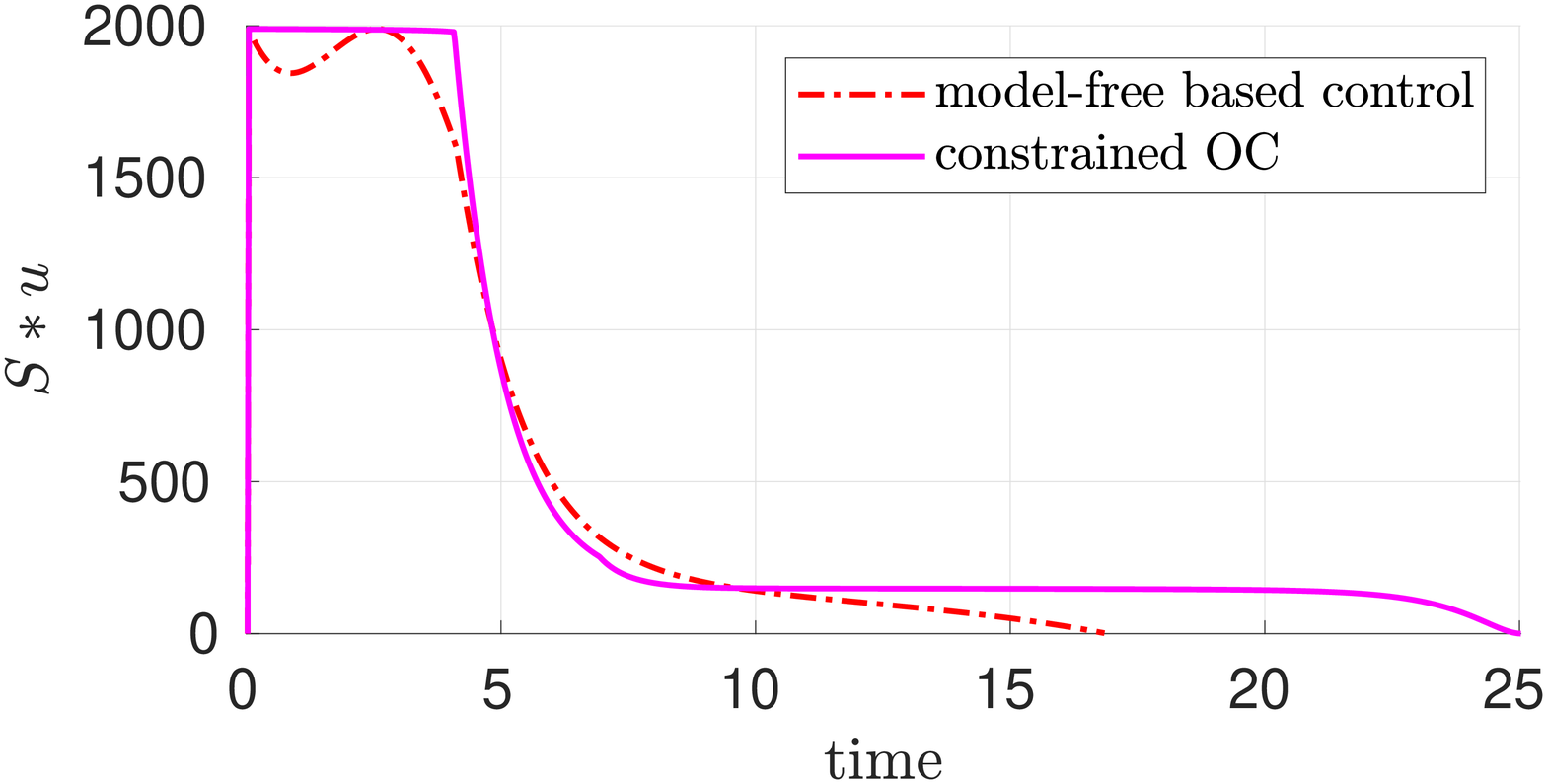}
\label{fg:opt_V:b}}
\caption{Evaluation of the constrained model-free based control: quadratic 
as first sequence including the constraint \eqref{eq:constraint:num} -- {\it case II}.}
\label{fg:opt_V}
\end{figure*}

In Table~\ref{tab:tab1}, we evaluate several cases with the cost criteria 
$J_{u+I}$,  $J_{I}$ and  $J_{u+I}^{Te}$, according to the final time 
of medication $T_e$. We compare the constrained and unconstrained classical optimal 
control problems and the unconstrained and constrained model-free problems 
with two types of configurations: \emph{slope} and \emph{quadratic} initial transient. 
The classical optimal control problem corresponds to the one performed 
in \cite{SilvaTorres:PrEP:DCDS:2018}.
\begin{table}[ht!]
\caption{Evaluation of the cost criteria.}
\label{tab:tab1}
\begin{tabular}{l|cccccccc} \hline
{Case}  & $T_e$ & $J_{u+I}$  &  $J_{I}$ &  $J_{u+I}^{Te}$ 
& $I(T_e)$ &  $\displaystyle{\max_{[0, t_{final}] } S u}$ & $u_{max}$ \\ \hline
{\small Unconstrained model-free} & 11.3 & $5.40 10^4$ & $5.40 10^4$ 
& $1.15 10^5$ & 31.12 & 3129 & 0.70 \\
{\small Constrained model-free -- slope (I)} 
& 19.0 & $6.45 10^4$ & $6.45 10^4$ & $2.39 10^5$ & 29.80 & 1990 & 0.80 \\
{\small Constrained model-free -- slope (II)} 
& 22.9 & $6.45 10^4$ & $6.45 10^4$ & $2.86 10^5$ & 28.25 & 2000 & 0.62 \\
{\small Constrained model-free -- quad. (I)} 
& 16.9 & $5.83 10^4$ & $5.83 10^4$ & $1.83 10^5$ & 29.12 & 1989 & 0.70 \\
{\small Constrained model-free -- quad. (II)} 
& 17.2 & $6.66 10^4$ & $6.66 10^4$ & $2.39 10^5$ & 32.29 & 1604 & 0.62  \\
{\small Unconstrained classical OC} & 25.0 
& $4.17 10^{4}$ & $4.17\, 10^{4}$ & $1.69 10^{5}$ & 21.95 & 9750 & 1 \\
{\small Constrained classical OC} 
& 25.0 & $6.14 10^{4}$ & $6.14\, 10^{4}$ & $2.72 10^{5}$ & 24.23 & 1989 & 1 \\ 
\hline
\end{tabular} 
\end{table}
Figures~\ref{fg:non_opt}--\ref{fg:opt_V} illustrate several scenarios: 
the constrained case is not satisfied under \emph{slope}, 
see Figure~\ref{fg:non_opt}; the constraint \eqref{eq:constraint:num} 
is satisfied under \emph{slope}, see Figures~\ref{fg:opt_II}--\ref{fg:opt_III}; 
the constrained cases are satisfied under \emph{quadratic function}, 
see Figures~\ref{fg:opt_IV}--\ref{fg:opt_V}.

From Figure~\ref{fg:non_opt:a}, we see that although the model-free based control 
and the classical optimal control approaches propose completely different control 
functions, the associated number of HIV infected individuals $I(t)$, $t \in [0, 25]$, 
are very similar. Interestingly, the control solution of the model-free based control 
is active in a much smaller interval of time that the one obtained with the classical 
optimal control approach: $T_e = 11.3$ \emph{versus} $T_e = 25$ (see first line of
Table~\ref{tab:tab1}). Analogous conclusions are taken 
from Figures~\ref{fg:opt_II:a}--\ref{fg:opt_V:a}.

It should be noted, however, that in Figure~\ref{fg:non_opt:b} the control obtained 
from the model-free based approach does not satisfy the mixed state-control 
constraint $S(t) u(t) \leq 2000$ for all $t$. In order to satisfy this constraint 
one must increase the time interval where the model-free control is active, 
see Figures~\ref{fg:opt_II:b}--\ref{fg:opt_V:b}. This depends on the configuration 
of the initial transient of the control (slope or quadratic), 
see Figures~\ref{fg:opt_II:a}--\ref{fg:opt_V:a} and Table~\ref{tab:tab1}.


\section{Discussion}
\label{sec:discussion}

The proposed control sequence can be considered as ``quasi-optimal'' in the sense 
that it does not obey to the Pontryagin maximum principle, so it is not an optimal control by definition, 
but it offers similar properties in terms of cost criteria minimization and reduction 
of the duration of treatment that is assured by our procedure.

The sequence is fully parametrized thanks to the initial transient coefficients $(L$ or $ Q)$ 
associated to $u_{max}$, including the $\mathcal{C}_{\pi}$-parameters for the decreasing transient, 
that must be adjusted according to the evolution of the infected state. In particular, 
the maximum value on the product $S u$ depends on the ``speed'' of the increasing transient 
as well as the final value $I(T_e)$, which depends on the ``speed'' of the increasing transient, 
the initial value $u_0$, and the feedback control that stabilizes $I$ 
through the decreasing transient. The transient slope plays a key role in the ``accuracy'' 
of the initial decrease of the infected state since a sufficient dose of the medication $u$ 
must be injected to the population in order to maintain the infected state to a lower level. 
Therefore, the maximum value of $u$ is a trade-off between the constraint $S u \leq \gamma_{max}$ 
to be satisfied and the duration of the treatment. Figure~\ref{fg:non_opt} illustrates 
a rather quick treatment, involving thus a fast initial transient but the constraint 
$S u \leq 2000$ is not satisfied; slower medical treatment for which the injection of the medication 
$u$ takes more time due to the constraint $S u \leq 2000$, can reduce the final asymptotic value $I(T_e)$ 
despite not necessarily fully reducing the cost criteria. The model-free based control aims to relax 
the treatment until $u = 0$ is reached. A first tuning has been made according to the gained experience 
and a more precise tuning can be performed using \cite{Porcelli}. The numerical evaluation 
of the cost criteria shows that our approach is globally better in terms of energy minimization. 
Moreover, the time-pondered criteria shows that the proposed control procedure is favorable 
to our approach taking into account the reduction of the duration of the treatment. These results 
illustrate afterwards that tuning the parameters of the proposed control sequence is a trade-off 
between considering minimizing the cost criteria, or minimizing the final value $I(T_e)$, 
depending of additional constraints.

We remark that the model-free based control could have been also used to drive the initial 
transient instead of the slope or the quadratic function, but the current algorithm offers 
slower performances at the very beginning to initiate the increase of the medication 
that prevents it to deal with, for example, the constraint $S u \leq 2000$.


\section{Conclusion and future work}
\label{sec:conclusion}

We have considered a SICAE epidemiological model for HIV/AIDS transmission,
proposing, for the first time in the literature, a model-free based 
approach to minimize the number of infected individuals. 
This approach consists in initializing PrEP medication, 
using a basic linear or quadratic function, and after that creating 
a direct feedback to control the decrease of infected individuals 
with respect to the considered measure of infected cases. 
Globally, the advantages of the proposed approach, 
when compared with the classical optimal control based on the 
Pontryagin Maximum Principle, is that it does not need any a priori knowledge 
of the model and a simple tuning of the proposed control 
sequence values allows good performances in terms of ``energy'' minimization 
as well as minimization of the medical treatment duration.
We concluded that our control strategy highlights interesting performances 
compared with the classical optimal control approach used 
in \cite{SilvaTorres:PrEP:DCDS:2018}. 

From a biological point of view, our application of the model-free based
control approach allows to propose new solutions for the implementation of PrEP in the prevention 
of HIV transmission, considering the constraints associated to the limitations on the availability 
of medicines for HIV and on number of individuals that the health systems have capacity to follow up
during their treatment. To the best of our knowledge, we were the first to apply the model-free 
based control approach in the context of epidemiology.

Future work may include replacing the slope or the quadratic 
initial transient by an optimal control; improvement of the proposed
model-free based control; implementation and comparison with 
the original Fliess--Join version of the model-free control,  
as it has been done, for example, for the glycemia control
\cite{FliessJoin_2021}. Stability issues regarding the closed loop
are very important and a promising LMI framework dedicated 
to study the stability of optimization algorithms is also of 
interest \cite{Fazlyab,sanzserna2021}.


\section*{Acknowledgments}

This work was partially supported by Portuguese funds 
through CIDMA, The Center for Research and Development 
in Mathematics and Applications of University of Aveiro, 
and the Portuguese Foundation for Science and Technology 
(FCT -- Funda\c{c}\~ao para a Ci\^encia e a Tecnologia), 
within project UIDB/04106/2020. Silva is also supported 
by the FCT Researcher Program CEEC Individual 2018 
with reference CEECIND/00564/2018.

The authors are grateful to Reviewers for several 
constructive suggestions and remarks.




\begin{thebibliography}{10}

\bibitem{Astrom} (MR3157723) [10.1016/j.automatica.2013.10.012]
\newblock K. J. Astr\"{o}m, P. R. Kumar, 
\newblock Control: A perspective,
\newblock \emph{Automatica},
\newblock \textbf{50} 1 (2014), 3--43.

\bibitem{Bara2} [10.1109/ECC.2016.7810602]
\newblock O. Bara, M. Fliess, C. Join, J. Day, S. Djouadi, 
\newblock Model-free immune therapy: A control approach to acute inflammation,
\newblock 15th European Control Conference (ECC'16), (2016). 

\bibitem{Bara} (MR3789959) [10.1016/j.jtbi.2018.04.003]
\newblock O. Bara, M. Fliess, C. Join, J. Day, S. M. Djouadi,
\newblock Toward a model-free feedback control synthesis for treating acute inflammation,
\newblock \emph{Journal of Theoretical Biology},
\newblock \textbf{448} (2018), 26--37.

\bibitem{cdc}
\newblock CDC,
\newblock Pre-Exposure Prophylaxis (PrEP),
\newblock Division of HIV/AIDS Prevention, 
National Center for HIV/AIDS, Viral Hepatitis, STD, and TB Prevention, 
Centers for Disease Control and Prevention,
\newblock \url{https://www.cdc.gov/hiv/risk/prep/index.html}, 
Page last reviewed: August 6, 2021.	 

\bibitem{AIDS:chronic:Lancet:2013} [10.1016/S0140-6736(13)61809-7]
\newblock S.~G.~Deeks, S.~R.~Lewin, D.~V.~Havlir,
\newblock The end of AIDS: HIV infection as a chronic disease,
\newblock \emph{The Lancet},
\newblock \textbf{382} 9903 (2013), 1525--1533.

\bibitem{D:S:T:2018} (MR3808514) [10.1016/j.aml.2018.05.005]
\newblock J. Djordjevic, C. J. Silva, D. F. M. Torres, 
\newblock A stochastic SICA epidemic model for HIV transmission, 
\newblock \emph{Appl. Math. Lett.} \textbf{84} (2018), 168--175.
\newblock {\tt arXiv:1805.01425}

\bibitem{Fazlyab} (MR3856216) [10.1137/17M1136845]
\newblock M. Fazlyab, A. Ribeiro, M. Morari, V.M. Preciado,
\newblock Analysis of Optimization Algorithms 
via Integral Quadratic Constraints: Nonstrongly Convex Problems,
\newblock \emph{SIAM Journal on Optimization},
\newblock \textbf{28} 3 (2018), 2654--2689.

\bibitem{Fliess2008a} 
\newblock M. Fliess, C. Join,
\newblock Commande sans mod\`{e}le et commande \`{a} mod\`{e}le restreint,
\newblock \emph{e-STA Sciences et Technologies de l'Automatique, 
SEE - Soci\'et\'e de l'Electricit\'e, de l'Electronique et des Technologies 
de l'Information et de la Communication},
\newblock  \textbf{5}(4) (2008), 1--23.

\bibitem{Fliess2009} [10.3182/20090706-3-FR-2004.00256]
\newblock M. Fliess, C. Join,
\newblock Model-free control and intelligent PID controllers: 
Towards a possible trivialization of nonlinear control?,
\newblock \emph{IFAC Proceedings Volumes}, 
\newblock  \textbf{42}(10) (2009), 1531--1550.

\bibitem{Fliess} (MR3172473) [10.1080/00207179.2013.810345]
\newblock M. Fliess, C. Join,
\newblock Model-free control,
\newblock \emph{International Journal of Control},
\newblock  \textbf{86} 12 (2013), 2228--2252.

\bibitem{FliessJoin_2021} [10.1002/rnc.5657]
\newblock M. Fliess, C. Join,
\newblock An alternative to PIs and PIDs: Intelligent proportional-derivative regulators, 
\newblock \emph{International Journal of Robust and Nonlinear Control}, in press.
\newblock DOI: 10.1002/rnc.5657

\bibitem{Hamiche} [10.1109/CoDIT.2019.8820297]
\newblock K. Hamiche, M. Fliess, C. Join, H. Aboua\"issa,
\newblock Bullwhip effect attenuation in supply chain management 
via control-theoretic tools and short-term forecasts: A preliminary study with an application to perishable inventories,
\newblock 6th International Conference on Control, Decision and Information Technologies (CoDIT),
\newblock 2019, 1492--1497.

\bibitem{Bernier} [10.1016/j.ifacol.2017.08.1167]
\newblock C. Join, J. Bernier, S. Mottelet, M. Fliess, 
S. Rechdaoui-Gu\'{e}rin, S. Azimi, V. Rocher, 
\newblock A simple and efficient feedback control strategy for wastewater denitrification, 
\newblock IFAC-PapersOnLine \textbf{50} (2017), no.~1, 7657--7662.

\bibitem{L:M:M:S:T:Y:2019} (MR3999700) [10.19139/soic.v7i3.834]
\newblock E. M. Lotfi, M. Mahrouf, M. Maziane, C. J. Silva, D. F. M. Torres and N. Yousfi, 
\newblock A minimal HIV-AIDS infection model with general incidence rate and application to Morocco data, 
\newblock \emph{Stat. Optim. Inf. Comput.} \textbf{7} (2019), no.~2, 588--603.
\newblock {\tt arXiv:1812.06965}

\bibitem{michel2018}
\newblock L. Michel,
\newblock \emph{A para-model agent for dynamical systems}, 
\newblock preprint {\tt arXiv:1202.4707}, (2018).

\bibitem{glycemia} [10.1109/TBME.2017.2698036]
\newblock T. MohammadRidha, M. A\"it-Ahmed, L. Chaillous, M. Krempf, I. Guilhem, J.-Y. Poirier, C. H. Moog,
\newblock Model Free iPID Control for Glycemia Regulation of Type-1 Diabetes, 
\newblock \emph{IEEE Transactions on Biomedical Engineering}, 
\newblock \textbf{65} 1 (2018), 199--206.

\bibitem{Ridha} [10.1137/1.9781611974072.9]
\newblock T. MohammadRidha, C. Moog, E. Delaleau, M. Fliess, C. Join, 
\newblock A variable reference trajectory for model-free glycemia regulation,
\newblock SIAM Conference on Control $\&$ its Applications (SIAM CT15), (2015).

\bibitem{Nascimento} [10.1016/j.arcontrol.2019.08.004]
\newblock T. P. Nascimento, M. Saska, 
\newblock Position and attitude control of multi-rotor aerial vehicles: A survey,
\newblock \emph{Annual Reviews in Control}, 
\newblock \textbf{48} (2019), 129--146.

\bibitem{Porcelli} [10.1145/3085592]
\newblock M. Porcelli, P.L. Toint,
\newblock BFO, A Trainable Derivative-Free Brute Force Optimizer 
for Nonlinear Bound-Constrained Optimization 
and Equilibrium Computations with Continuous and Discrete Variables,
\newblock \emph{ACM Transactions on Mathematical Software}, 
\newblock \textbf{44} 1, no.~6 (2017), 1--25.

\bibitem{sanzserna2021} (MR4267491) [10.1137/20M1364138]
\newblock J. M. Sanz-Serna, K. C. Zygalakis,
\newblock The connections between Lyapunov functions 
for some optimization algorithms and differential equations, 
\newblock \emph{SIAM J. Numer. Anal.},
\newblock \textbf{59} (3) (2021), 1542--1565.

\bibitem{SilvaTorres:TBHIV:2015} (MR3392642) [10.3934/dcds.2015.35.4639]
\newblock C. J. Silva, D. F. M. Torres,
\newblock A TB-HIV/AIDS coinfection model and optimal control treatment,
\newblock \emph{Discrete Contin. Dyn. Syst.},
\newblock \textbf{35} no.~9 (2015), 4639--4663.
\newblock {\tt arXiv:1501.03322}

\bibitem{SilvaTorres:EcoComplexity} [10.1016/j.ecocom.2016.12.001]
\newblock C. J. Silva, D. F. M. Torres,
\newblock A SICA compartmental model in epidemiology with application to HIV/AIDS in Cape Verde, 
\newblock \emph{Ecological Complexity},
\newblock \textbf{30} (2017), 70--75.
\newblock {\tt arXiv:1612.00732}

\bibitem{SilvaTorres:PrEP:DCDS:2018} (MR3714435) [10.3934/dcdss.2018008]
\newblock C. J. Silva, D. F. M. Torres,
\newblock Modeling and optimal control of HIV/AIDS prevention through PrEP,
\newblock \emph{Discrete Contin. Dyn. Syst. Ser. S},
\newblock \textbf{11} no. 1 (2018), 119--141.
\newblock {\tt arXiv:1703.06446}

\bibitem{S:T:2019} (MR3980195) [10.1016/j.matcom.2019.03.016]
\newblock C. J. Silva and D. F. M. Torres, 
\newblock Stability of a fractional HIV/AIDS model, 
\newblock \emph{Math. Comput. Simul.} \textbf{164} (2019), 180--190. 
\newblock {\tt arXiv:1903.02534}

\bibitem{MyID:455} (4175342) [10.1007/978-3-030-49896-2_6]
\newblock C. J. Silva and D. F. M. Torres,
\newblock On SICA models for HIV transmission.
\newblock In: \emph{Mathematical Modelling and Analysis of Infectious Diseases, 
Studies in Systems, Decision and Control} \textbf{302} (2020), 155--179.
\newblock {\tt arXiv:2004.11903}

\bibitem{Tebbani}
\newblock S. Tebbani, M. Titica, C. Join, M. Fliess, D. Dumur, 
\newblock \emph{Model-based versus model-free control designs 
for improving microalgae growth in a closed photobioreactor: Some preliminary comparaisons},
\newblock The 24th Mediterranean Conference on Control and Automation (MED'16), IEEE (2016), 683--688.

\bibitem{V:T:2021} (MR4276109) [10.3934/mbe.2021231]
\newblock S. Vaz\ and\ D. F. M. Torres, 
\newblock A dynamically-consistent nonstandard finite difference scheme for the SICA model, 
\newblock \emph{Math. Biosci. Eng.} \textbf{18} (2021), no.~4, 4552--4571. 
\newblock {\tt arXiv:2105.10826}

\bibitem{MR4261252} (MR4261252) [10.1186/s13662-021-03392-y]
\newblock X. Wang, C. Wang\ and\ K. Wang, 
\newblock Extinction and persistence of a stochastic SICA epidemic model 
with standard incidence rate for HIV transmission, 
\newblock {\emph Adv. Difference Equ.} \textbf{2021}, Paper No.~260, 17~pp. 

\bibitem{art:viral:load} [10.1016/S0140-6736(08)61115-0]
\newblock D. P. Wilson, M. G. Law, A. E. Grulich, D. A. Cooper, J. M. Kaldor,
\newblock Relation between HIV viral load and infectiousness: 
A model-based analysis,
\newblock \emph{The Lancet},
\newblock \textbf{372} no.~9635 (2008), 314--320.

\bibitem{MR4279814} (MR4279814) [10.1002/mma.7292]
\newblock M. Zaka Ullah\ and\ D. Baleanu, 
\newblock A new fractional SICA model and numerical method for the transmission of HIV/AIDS, 
\newblock \emph{Math. Methods Appl. Sci.} \textbf{44} (2021), no.~11, 8648--8659. 

\end{thebibliography}
\end{document}